\documentclass{article}
\usepackage{amssymb}
\usepackage{a4}
\usepackage{epsfig}
\usepackage{latexsym}
\usepackage{pictex}
\usepackage{graphicx}
\usepackage{amsmath}
\usepackage{amsfonts}
\begin{document}
%%%%%%%%%%%%%%%%%%%%%%%%%%%%%%%%%%%%%%%%%%%%%%%%%%%%%%%%%%%%%%%%%%%%%%
%	spaces for your own definitions follows
%%%%%%%%%%%%%%%%%%%%%%%%%%%%%%%%%%%%%%%%%%%%%%%%%%%%%%%%%%%%%%%%%%%%%%

\providecommand{\BBb}[1]{{\mathbb{#1}}}
\providecommand{\cal}[1]{{\mathcal{#1}}}   

\newcommand{\Dm}{\BBb D}
\newcommand{\fracc}[2]{{\textstyle\frac{#1}{\raise 1pt\hbox{$\scriptstyle #2$}}}}
\newcommand{\fracnp}{\fracc np}

\newcommand{\g}{\gamma_0}
\newcommand{\lap}{\operatorname{\Delta}}

\newcommand{\mlap}{-\!\operatorname{\Delta}}
\newcommand{\N}{\BBb N}
\newcommand{\ran}{\operatorname{ran}}
\newcommand{\R}{{\BBb R}}
\newcommand{\Rn}{{\BBb R}^{n}}
\newcommand{\singsupp}{\operatorname{sing\,supp}}
\newcommand{\supp}{\operatorname{supp}}

\newcommand{\dom}{\operatorname{dom}}

\newtheorem{thm}{Theorem}[section]
\newtheorem{cor}[thm]{Corollary}
\newtheorem{rem}{Remark}[section]

%%%%%%%%%%%%%%%%%%%%%%%%%%%%%%%%%%%%%%%%%%%%%%%%%%%%%%%%%%%%%%%%%%%%%%
%	Please fill in your details below,
%	address at the end of the file
%%%%%%%%%%%%%%%%%%%%%%%%%%%%%%%%%%%%%%%%%%%%%%%%%%%%%%%%%%%%%%%%%%%%%%

\author{Jon Johnsen}

\title{Regularity results and parametrices of semi-linear 
boundary problems of product type.}

\date{\em Dedicated to Prof. Hans Triebel on the occasion of his 65th
birthday}

%\date{}

\maketitle

%%%%%%%%%%%%%%%%%%%%%%%%%%%%%%%%%%%%%%%%%%%%%%%%%%%%%%%%%%%%%%%%%%%%%%%
%	You can insert an abstract here if you want. 
%%%%%%%%%%%%%%%%%%%%%%%%%%%%%%%%%%%%%%%%%%%%%%%%%%%%%%%%%%%%%%%%%%%%%%%
%\begin{abstract}

%\end{abstract}
%%%%%%%%%%%%%%%%%%%%%%%%%%%%%%%%%%%%%%%%%%%%%%%%%%%%%%%%%%%%%%%%%%%%%%%
%	Please insert the article body now
%%%%%%%%%%%%%%%%%%%%%%%%%%%%%%%%%%%%%%%%%%%%%%%%%%%%%%%%%%%%%%%%%%%%%%%
\section{Introduction}
  \label{intr-sect}
This study focuses on semi-linear problems of the form
\begin{equation}
\begin{aligned}
    Au+ N(u)&=f \quad\text{in}\quad \Omega 
 \\ 
    Tu&=\varphi\quad\text{on}\quad \Gamma:=\partial\Omega.
\end{aligned}
  \label{mpb-eq}  
\end{equation}
Here $(f,\varphi)$ are the given data, and $u$ the unknown. Problem
\eqref{mpb-eq} should 
be elliptic in some bounded, $C^\infty$-smooth region $\Omega\subset\Rn$;
that is $A$ should be a 
linear differential operator in $\Omega$ while $T$ should be a
trace operator such that the system $\{A,T \}$ is elliptic in
$\Omega$. More generally, $A$ could be suitably ``pseudo-differential'' as
long as $\{A,T \}$ is injectively elliptic in the Boutet de~Monvel
calculus of boundary problems.

$N(u)$ stands for a non-linearity which combines $u(x)$ and its
derivatives $D^\alpha u$ in a polynomial way, roughly speaking.
 
The main point is the following frequently asked question:  given a solution
$u$, does the presence of $N(u)$ influence the regularity of $u$?

This problem can of course be phrased in various frameworks:
to measure regularity, the Besov and Triebel--Lizorkin spaces
$B^{s}_{p,q}$ and $F^{s}_{p,q}$ could be adopted for $s\in \R$ and $p$,
$q\in\,]0,\infty]$ (though with $p<\infty$ for $F^{s}_{p,q}$). 
But to simplify matters\,---\,and indeed to fix ideas\,---\,this
survey deals with the Sobolev, or Bessel potential, spaces
$H^s_p(\overline{\Omega})$, 
where $s\in\R$ and $1<p<\infty$.
Now the solution may be known to exist in some a priori space, denoted
$H^{s_0}_{p_0}$ throughout, while data are given in other spaces having
some integral exponent $r\in\,]1,\infty[\,$. The case with $r\neq p_0$
requires extra efforts, and the present paper deals with a flexible way of
handling this.

\bigskip

The word ``semi-linear'' is often taken as an indication that solutions of
problems like \eqref{mpb-eq} will have practically the same regularity as in
the linear case, ie as when $N\equiv0$. However, when $r\neq p_0$ is
allowed, it is more demanding to describe for which a priori spaces and data
spaces this property of semi-linearity holds.

A classical way to obtain such conclusions is to improve
the knowledge of $u$ in finitely many steps (ie a boot-strap method). 
But one faces rather pains-taking difficulties when this method is applied
to cases in which the a priori space for $u$
is not ``close enough'' to the solution space associated to the data $f$ and
$\varphi$ in the linear theory.
(Such phenomena have been described in \cite{JJ93,JJ95stjm} and in a joint
work with T.~Runst \cite{JoRu97}).

However, in a recent article \cite{JJ01} a different technique was 
worked out\,---\,it requires rather weaker assumptions than
boot-strap methods do, it has cleaner proofs and in particular 
it also avoids the technicalities mentioned above.
In short this approach is a much more flexible tool.
It was exemplified for elliptic problems in full generality
in \cite{JJ01}, where the crucial point was 
a specific parametrix formula for the
non-linear problem \eqref{mpb-eq}; this formula is useful because one can
read off a given solution's regularity directly. 

The purpose of the present paper is to give a concise account of the
resulting technique and to present how the parametrices
straightforwardly give regularity improvements.

To give a very brief account of the outcome of the study (with examples to
follow further below),  
it is useful to introduce three \emph{parameter domains}:
\begin{equation}
  \Dm({\cal A}),\quad \Dm(N),\quad \Dm(L_u).
  \label{Dm-eq}
\end{equation}
Here $\Dm({\cal A})$ consists of all the (pairs of) parameters $(s,p)$ for
which the matrix-formed operator ${\cal A}:=
\left(\begin{smallmatrix}A\\T\end{smallmatrix}\right)$ is defined on the
space $H^s_p$.  
This takes into account the class of $T$ (and of $A$ in the
pseudo-differential case).

Similarly $\Dm(N)$
contains all the $(s,p)$ for which $N$ is defined on $H^s_p$
and has order \emph{less} than that of $A$ on
this space. Finally, and most importantly, for any $(s_0,p_0)\in\Dm(N)$ and
any given $u$ in $H^{s_0}_{p_0}$, there should exist some linear but
possibly $u$-dependent operator $L_u$ such that
\begin{equation}
  N(u)=-L_u(u).
  \label{NLu-eq}
\end{equation} 
When the operator $L_u$ is studied in its own right on the scale
$H^{s}_{p}$, then $\Dm(L_u)$ contains the $(s,p)$ for which $L_u$ is
defined and has lower order than $A$. In addition it is necessary to require
of $N$ and $L_u$
that $\Dm(N)\subset \Dm(L_u)$.

In practice $\Dm(L_u)$ is much larger than $\Dm(N)$, and the more regular
$u$ is known a priori to be, the larger $\Dm(L_u)$ will be. (While $\Dm(N)$
is the same independently of any given solution $u$.) This leads to a main
feature:

On the one hand, using a boot-strap method it turns out that one can work
inside the domain $D(A)\cap \Dm(N)$; which is logical because $N$ would
lose ``more derivatives'' on spaces outside $\Dm(N)$.
On the other hand, the present parametrix methods work well on the larger
set
\begin{equation}
  \Dm_u:=\Dm({\cal A})\cap\Dm(L_u).
\end{equation}
For this reason a given solution may be treated under much weaker initial
assumptions on the data $(f,\varphi)$. Indeed, the regularity of $u$
is read off from the following \emph{parametrix formula} (derived in
\cite{JJ01})
\begin{equation}
  u=P^{(N)}_u(Rf+K\varphi)+{\cal R}u+(RL_u)^{N}u.
  \label{para-eq}
\end{equation}
Here $\left(\begin{smallmatrix}R&K\end{smallmatrix}\right)$ denotes a 
left-parametrix of $\left(\begin{smallmatrix}A\\T\end{smallmatrix}\right)$,
with associated smoothing operator ${\cal R}$; that is 
$RA+KT=I-{\cal R}$. The parametrix
$P^{(N)}_u$ is a finite Neumann series with the linear operator $RL_u$
as `quotient'.
Consequently, with known mapping properties of $R$, $K$ and $L_u$, the above
formula \eqref{para-eq} shows directly how the regularity of $u$ is
determined by the data together with the a priori regularity of $u$ itself
(the latter enters the term $(RL_u)^N$).

Below this is explained in detail by means of an example.

\section{The Framework}
  \label{frmw-sect}
As a another simplification we may consider the following model problem,
which is rich enough to illustrate the points. Here and below
$\gamma_0$ denotes the standard trace (restriction) on $\Gamma$:
\begin{equation}
\begin{aligned}
  \mlap u +u\partial_{x_1}u&=f  \quad\text{in}\quad\Omega
\\
  \g u&=\varphi \quad\text{on}\quad\Gamma.
\end{aligned}
  \label{mdp-eq}  
\end{equation}
The discussion of \eqref{mdp-eq} will be carried out under the hypothesis
that a solution $u$ is given for some specific data $(f,\varphi)$ fulfilling
\begin{align}
  u&\in H^{s_0}_{p_0}(\overline{\Omega})
  \label{uass-eq}  \\
  f&\in H^{t-2}_r(\overline{\Omega}),\qquad \varphi\in
B^{t-1/r}_{r,r}(\Gamma).
  \label{dass-eq}
\end{align}
In general, the space $H^s_p(\overline{\Omega})$ is defined by restriction to
$\Omega$ and $B^s_{p,p}(\Gamma)$ is defined via local coordinates on
$\Gamma$.
 
With this set-up, the theme is whether $u$ belongs to the space
$H^{t}_r(\overline{\Omega})$ too. (There are of course necessary
conditions for this to be true, eg $s_0>1/p_0$ must hold for the boundary
condition to make sense. It is tempting to require
in analogy that $t>1/r$, but it is a point that weaker assumptions will
suffice; hence this discussion is postponed a little.) 
Using boot-strap arguments in treating this, difficulties occur as mentioned
in the introduction. Indeed,
for cases with, say $p_0$ and $r$ close to $1$ and $\infty$ respectively,
and small values of $s_0>t$, boot-strapping is possible, but
careful arguments based on special estimates of $u\partial_1 u$ are needed
to avoid auxiliary spaces on which $\gamma_0 u$ is undefined; cf
\cite{JJ95stjm}. 

With a more direct approach, 
the aim is to ``invert'' \eqref{mdp-eq} by means of the formula
\begin{equation}
  u=P^{(N)}_u(R_Df+K_D\varphi)+{\cal R}u+(R_DL_u)^{N}u.
  \label{para-id}
\end{equation}
To explain the various quantities in \eqref{para-id}, it is first noted that
the linear problem corresponding to \eqref{mdp-eq} is considered as
an equation for the elliptic Green operator (when $s>1/p$, $1<p<\infty$),
\begin{equation}
  {\cal A}= \begin{pmatrix}\mlap\\ \g\end{pmatrix}\colon 
  H^s_p(\overline{\Omega})\to 
  \begin{gathered}
    H^{s-2}_p(\overline{\Omega})
  \\
     \oplus 
  \\
    B^{s-1/p}_{p,p}(\Gamma).
  \end{gathered}
  \label{Dir-eq}
\end{equation}
Then $\left(\begin{smallmatrix} R_D & K_D\end{smallmatrix}\right)$
is taken as a parametrix (belonging to the Boutet de~Monvel calculus), 
ie it is continuous in the opposite direction in \eqref{Dir-eq} and 
\begin{equation}
   \begin{pmatrix} R_D& K_D\end{pmatrix}
      \begin{pmatrix}\mlap\\ \g\end{pmatrix} =I-{\cal R};
\end{equation}
here the range ${\cal R}(H^s_p)\subset C^\infty(\overline{\Omega})$ for all
$s>1/p$. (In fact ${\cal R}=0$ is possible for this Dirichl\'et problem, but
it is retained here to make it clear that also in general its presence is
harmless.) 

The second ingredient in \eqref{para-id} is a decomposition of the
non-linear term as 
\begin{equation}
  u\partial_{x_1}u =-L_u(u),\qquad\text{$L_u$ linear}.
  \label{Lu-eq}
\end{equation}
More precisely, it is necessary to ensure that $L_u$ has certain mapping
properties, hence it is defined by means of a universal extension operator
from $\Omega$ to $\Rn$, say $\ell_\Omega$, to be
\begin{equation}
  -L_u(v)=\pi_1(\ell_\Omega u, \ell_\Omega\partial_1v)+
  \pi_2(\ell_\Omega u, \ell_\Omega\partial_1v)+
  \pi_3(\ell_\Omega v, \ell_\Omega\partial_1u).
  \label{Luv-eq}
\end{equation}
Here the $\pi_j(\cdot,\cdot)$ are para-multiplication operators defined on
$\Rn$ (so that restriction to $\Omega$ 
of each term on the right hand side of \eqref{Luv-eq} is understood). 
These are introduced using a Littlewood--Paley partition
$1=\sum_{j=0}^\infty \Phi_j(\xi)$ with smooth functions $\Phi_j$ supported
at $\{\,2^{j-1}\le|\xi|\le2^{j+1} \,\}$ for $j>0$; then
\begin{equation}
  \pi_1(g,h)=\sum_{j=2}^\infty (\Phi_0(D)+\dots+\Phi_{j-2}(D))g\cdot
   \Phi_j(D)h,
\end{equation}
and $\pi_3(g,h):=\pi_1(h,g)$ whilst $\pi_2(g,h)$ gives the remainder in the
formal decomposition of $g\cdot h$.

Using this, the parametrices $P^{(N)}_u$ of the non-linear problem
\eqref{mdp-eq} are now finally introduced as
\begin{equation}
  P^{(N)}_u = I +R_DL_u +\dots+(R_DL_u)^{(N-1)}, \qquad N\in\N.
  \label{PNu-eq}
\end{equation}
They clearly depend on the given solution $u$, and since
both $R_D$ and $L_u$ are linear, so are the $P^{(N)}_u$ on every
$H^s_p(\overline{\Omega})$ with $(s,p)\in\Dm({\cal A})\cap\Dm(L_u)$; 
cf~\eqref{DmA-eq} ff.

For the above model problem, the parameter domains from the introduction
are,
when $N(u)=u\partial_1 u$ and $t_+=\max(0,t)$ denotes the positive part, 
\begin{align}
  \Dm({\cal A})&=\{\,(s,p)\mid s>1/p \,\}
  \label{DmA-eq}  \\
  \Dm(N)&=\{\,(s,p)\mid s>\tfrac{1}{2}
     +(\fracnp-\tfrac{n}{2})_+,\quad s>\fracnp-1 \,\}
  \label{DmN-eq} \\
  \Dm(L_u)&=\{\,(s,p)\mid s+s_0>1
     +(\fracc n{p_0}+\fracnp-n)_+ \,\}
  \label{DmL-eq} 
\end{align}
The two first restrictions make the trace $\gamma_0u$ and the product
$u\cdot \partial_1 u$ well defined on $H^s_p(\overline{\Omega})$, whilst the
second condition in \eqref{DmN-eq} implies that for some $s'>s-2$ it holds
that $u\partial_1u\in H^{s'}_p(\overline{\Omega})$ for every $u\in
H^s_p(\overline{\Omega})$. 

More noteworthy is it that the requirement in \eqref{DmL-eq} is effectively
weaker the better the a priori knowledge of $u$ is: for higher values of
$s_0$ or larger values of $p_0$, more pairs $(s,p)$ fulfil the inequality.

A closer analysis shows that
$L_u$ has order $\omega$ in the sense that
\begin{gather}
    L_u\colon H^s_p(\overline{\Omega})\to H^{s-\omega}_p(\overline{\Omega})
  \quad\text{for all $(s,p)\in\Dm(L_u)$,}
  \label{Lom-eq}  \\
  \omega=1+(\fracc n{p_0}-s_0)_++\varepsilon,\qquad \varepsilon\ge0.
  \label{om-eq}
\end{gather}
Here $\varepsilon>0$ is only necessary for $s_0=\fracc n{p_0}$.  If
one removes $\ell_\Omega$ and the restriction to $\Omega$ from $L_u$, then
the resulting operator is in the H{\"o}rmander class
$\operatorname{OPS}^\omega_{1,1}(\Rn\times\Rn)$, leading to an analogous
continuity property.

It is important to observe that with $u\in H^{s_0}_{p_0}(\overline{\Omega})$
for some $(s_0,p_0)$ in $\Dm(N)$, it follows from \eqref{DmN-eq} that the
order of $L_u$ satisfies $\omega<2$ in \eqref{Lom-eq}. In other words, $L_u$
loses fewer derivatives than $\lap$. It is a peculiar fact that $L_u$, once
$u$ is chosen, actually has constant order on all spaces regardless of
whether they have parameter inside or outside $\Dm(N)$ (by comparison, the
boundary of $\Dm(N)$ contains points where $N$ attains the order $2$).

Using the above continuity results for $R_D$, $K_D$ and $L_u$, one can now
show that the parametrix formula holds and derive the regularity results.

\begin{rem}
The operator $L_u$ in \eqref{Luv-eq} differs from the linearisations in
J.~M.~Bony's work \cite{Bon} because the $\pi_2$-term is a part of the
operator instead of being treated as a negligible error term. In the present
context this has to be so, for it does occur that this term has a
non-negligible regularity and in addition it would not be natural to violate
the identity $L_uu=u\partial_1 u$. For this reason it is suggested that
one could call $L_u$ the \emph{full} paralinearisation of $u\partial_1 u$.

Moreover, the perhaps more natural linearisations  $v\mapsto
u\partial_1 v$ and the differential $u\partial_1 v+v\partial_1 u$ do not work
in this context, because they are not moderate in the terminology of
\cite{JJ01}. Indeed, on $H^s_p$ they have order equal
to $s-1-s_0$ for large $s$, and this has no upper limits for $s\to\infty$;
unlike $L_u$ that has constant order with respect to $s$ as observed above.
\end{rem}

\begin{rem}
About the above results it should be mentioned that the properties of linear
elliptic problems were deduced for the $H^s_p$-scale in full generality by
G.~Grubb \cite{G3}, who extended the Boutet de~Monvel calculus to these
spaces (and to the classical Besov spaces). In particular this implies 
\eqref{Dir-eq} and the statements following it. (For the $B^{s}_{p,q}$ and
$F^{s}_{p,q}$ scales there is a similar extension of the calculus in
\cite{JJ96ell}, which applies to the present problems in the same way.)
For introductions to the calculus the reader may consult \cite{G97,G2}.

In the definition of $L_u$, the universal extension operator was constructed
by V.~Rychkov \cite{Ry99,Ry99BPT}. 
He showed that $\ell_\Omega$ can be taken such
that for all $s$ and $p$ it is continuous 
\begin{equation}
  \ell_\Omega\colon H^s_p(\overline{\Omega})\to H^s_p(\Rn)
  \label{lOm-eq}
\end{equation}
and that $r_\Omega\ell_\Omega=I$ holds on $H^s_p(\overline{\Omega})$
(in fact it was carried out for the Besov and Triebel--Lizorkin scales).

The para-multiplication operators $\pi_j(\cdot,\cdot)$ in \eqref{Luv-eq}
follow M.~Yamazaki \cite{Y1} in the notation and the definition. To prove 
\eqref{Lom-eq}--\eqref{om-eq} it suffices to combine \eqref{lOm-eq} with
standard estimates of the $\pi_j(\cdot,\cdot)$; these are essentially found
in \cite{Y1}, but for a proof the reader may consult \cite{JJ01}, which
presents a general study of non-linear operators of product type
(encompassing sums of terms $D^\alpha u\cdot D^\beta u$ and more general
expressions). For a full set of estimates of para-multiplication operators
proved directly (without the somewhat heavier paradifferential techniques in
\cite{Y1}), the reader may eg consult \cite[Th~5.1]{JJ94mlt}.

\end{rem}

\section{Results involving parametrices}
  \label{res-sect}

We shall now proceed to state the results for the model problem, that one
obtains from the parametrices. Recall that ${\cal A}=
\left(\begin{smallmatrix} \mlap\\ \gamma_0\end{smallmatrix}\right)$ and
$N(u)=u\partial_1 u$. Furthermore $\Dm({\cal A})$, $\Dm(N)$ and $\Dm(L_u)$
are given as in \eqref{DmA-eq} ff, so that they fulfil the conditions
described after \eqref{Dm-eq}.

The first point is to establish that the parametrix formula \eqref{para-id}
really holds, and to give basic properties of the entering operators.

\begin{thm}
  \label{para-thm}
Assume that \eqref{mdp-eq} and \eqref{uass-eq}--\eqref{dass-eq} hold for
parameters fulfilling 
\begin{equation}
  \begin{aligned}
  (s_0,p_0)&\in \Dm({\cal A})\cap 
\Dm(N)
  \\
(t,r)&\in \Dm({\cal A})\cap \Dm(L_u)=:\Dm_u.
\end{aligned}
  \label{pass-eq}
\end{equation}
Then \eqref{para-id} holds, and for $P^{(N)}_u$ as in \eqref{PNu-eq},
\begin{alignat}{2}
  \forall N\in\N_0,\ \forall (s,p)&\in \Dm_u\colon &
  P^{(N)}_u&\colon H^s_p(\overline{\Omega})
  \to H^s_p(\overline{\Omega}),
  \label{allN-eq}
  \\
  \exists N&\in\N_0\colon& \quad(R_DL_u)^{N}&\colon 
   H^{s_0}_{p_0}(\overline{\Omega})\to H^t_r(\overline{\Omega}).
  \label{exN-eq}
\end{alignat}
In these lines the arrows stand for continuous, linear maps.
\end{thm} 
Actually \eqref{exN-eq} holds for all sufficiently large values of $N$, but
usually it is enough have a single such $N$.

The above Theorem~\ref{para-thm} is a special case of an abstract result
proved in \cite[Th~2.2]{JJ01}. The proof is not difficult in itself; it is
formulated for a general situation specified by some lengthy, but
essentially rather mild conditions labelled (I)--(V) in
\cite[Sect.~2]{JJ01}, and that these are fulfilled for the model problem
considered in this paper is a consequence of the above Section~\ref{frmw-sect}.

Now one immediately gets

\begin{cor}
  \label{reg1-cor}
If \eqref{mdp-eq}, \eqref{uass-eq}--\eqref{dass-eq} and \eqref{pass-eq} all
hold, then $u$ is also an element of $H^t_r(\overline{\Omega})$.
\end{cor}

It is a mjor point of the paper that, to prove this, one may take $N$ as in
\eqref{exN-eq}; then the properties in \eqref{exN-eq}, 
\eqref{allN-eq}, \eqref{Lom-eq} together with formula \eqref{para-id} show
that $u\in H^t_r(\overline{\Omega})$.

It deserves to be emphasised that $(t,r)$ is assumed to lie in $\Dm(L_u)$
but not necessarily in the smaller parameter domain for the non-linear term,
$\Dm(N)$. For this reason it is possible to conclude  that given solutions
may belong to spaces that are beyond the reach of the boot-strap method.

\section{Final remarks}
This paper focuses on semi-linear elliptic boundary problems, and even
specialises to a simple model problem in order not to burden the exposition;
the possible extensions are many, but the reader should get a good
impression of the possibilities from the above. 
Within the framework of elliptic problems
some of the generalisations are indicated after \eqref{mpb-eq}, but is is
also possible to include semi-linear elliptic systems like the stationary
Navier--Stokes equation or von Karman's equations. This requires an extended
notion of product type non-linearities, defined on sections of vector
bundles. The reader may consult \cite{JJ01} for this. 

The general study in \cite[Sect.~2]{JJ01} also allow some applications to
parabolic initial-boundary problems with non-linearities of product
type. However, when such problems are non-homogeneous the compatibility
conditions on the data set severe restrictions to how much the regularity
can be improved; but even so it should be possible to work out some results
in this area.

Concerning the tools, it is on the one hand clear that it is a fine theory
of linear elliptic problems that enter, namely that of the Boutet de~Monvel
calculus. On the other hand, the treatment of the non-linear terms is based
on para-multiplication operators on $\Rn$. This technique was essentially
introduced (independently) by J.~Peetre and H.~Triebel around 1976-77
\cite{Pee,T-pmlt,T0} in order to analyse the pointwise product. One way to
sum up the present paper could be to say that para-multiplication also
may enter in a crucial way in treatments of certain non-linear 
perturbations of elliptic boundary problems. 

%%%%%%%%%%%%%%%%%%%%%%%%%%%%%%%%%%%%%%%%%%%%%%%%%%%%%%%%%%%%%%%%%%%%%%%
%	Please insert your list of references below
%%%%%%%%%%%%%%%%%%%%%%%%%%%%%%%%%%%%%%%%%%%%%%%%%%%%%%%%%%%%%%%%%%%%%%%

\vspace*{\fill}

%%%%%%%%%%%%%%%%%%%%%%%%%%%%%%%%%%%%%%%%%%%%%%%%%%%%%%%%%%%%%%%%%%%%%%%
%	Finally, your address, please
%%%%%%%%%%%%%%%%%%%%%%%%%%%%%%%%%%%%%%%%%%%%%%%%%%%%%%%%%%%%%%%%%%%%%%%

\noindent
Jon Johnsen \\
Department of Mathematical Sciences\\
Aalborg University\\
Fredrik Bajers Vej 7G\\
DK--9220 Aalborg East\\
Denmark\\

\noindent
jjohnsen@math.auc.dk

\end{document}